\documentclass[a4paper,draft]{amsart}

\usepackage{amsmath}
\usepackage{amssymb,amsthm,amscd}

\newtheorem{theorem}{Theorem}[section]
\newtheorem{lemma}[theorem]{Lemma}
\newtheorem{theo}[theorem]{Theorem}
\newtheorem{proposition}[theorem]{Proposition}

\theoremstyle{definition}

\newtheorem{example}[theorem]{Example}
\theoremstyle{remark}

\theoremstyle{remark}

\numberwithin{equation}{section}

\newcounter{count}
\newcounter{countroman}

\renewcommand{\le}{\leqslant}
\renewcommand{\ge}{\geqslant}
\def\boplus{{{\textstyle\bigoplus}}}

\def\coloneq{\mathrel{\mathop:}=}

\def\T{{\rm T}}

\def\eps{{\varepsilon}}
\def\bangle#1{\langle #1 \rangle}

\def\K{{\mathrm K}}
\def\z{\mathfrak{z}}
\def\LL{{\mathrm L}}
\def\ZZ{{\mathbb Z}}
\def\NN{{\mathbb N}}

\def\id{{\rm id}}
\def\Aut{\operatorname{Aut}}
\def\End{\operatorname{End}}

\def\Sym{{\rm Sym}}
\def\vect{{\rm vect}}

\def\Mat{{\rm Mat}}

\def\chr{{\rm char\:}}

\def\tr{{\rm tr \,}}

\def\embed{\hookrightarrow}

\def\Aut{{\rm Aut}}

\begin{document}
\sloppy

\author[A. Perepechko]%
{Alexander Perepechko}
\title[Affine monoids as monoids of endomorphisms]
{Affine algebraic monoids as endomorphisms' monoids of finite-dimensional algebras}
\address{Department of Higher Algebra, Faculty of Mechanics and Mathematics, Moscow State University, Leninskie Gory, Moscow, 119991, Russia}
\email{shuubun@gmail.com}

\begin{abstract}
We prove that any affine algebraic monoid can be obtained as the endomorphisms' monoid of a finite-dimensional (nonassociative) algebra.
\end{abstract}

\maketitle

\section{Introduction}
 Let $\K$ be an algebraically closed field of arbitrary characteristic. Recall that an \emph{affine algebraic semigroup} is an affine variety $M$ over $\K$ with an associative product $\mu\colon M\times M\to M$, which is a morphism of algebraic varieties. Denote an element $\mu(a,b)$ by $ab$. A semigroup is called a \emph{monoid} if it contains an identity element $e\in M$ such that $em=me=m$ for any $m\in M$. An element $0\in M$ is called \emph{zero} if $0m=m0=0$ for any $m\in M$. Obviously, a monoid cannot contain more than one zero. It is well known that every affine algebraic monoid is isomorphic to a Zariski closed submonoid of the monoid $\LL(V)$ of all linear operators on 
  some finite-dimensional vector space $V$, e.g. see \cite[Theorem 3.8]{R} or \cite[Lemma 1.11]{A}. A systematic account of the theory of affine algebraic monoids is given in \cite{P} and \cite{R}.
The classification of irreducible affine monoids, whose unit group is reductive, is obtained in \cite{Rit} and \cite{V}. 

Let $A$ be a finite-dimensional algebra over the field $\K$, i.e. a finite-dimensional vector space $A$ with a bilinear map $\alpha\colon A\times A\to A$. Note that the associativity or commutativity of the map $\alpha$ is not assumed. It is convenient to denote by $\vect(A)$ the underlying vector space of an algebra $A$. By an ideal of an algebra $A$ we mean a two-sided ideal. An algebra $A$ is called \emph{simple} if it does not contain proper ideals. The set of all endomorphisms of $A$,
$$
\End(A)\coloneq \{\phi\in\LL(\vect(A))\;|\;\alpha(\phi(a),\phi(b))=\phi(\alpha(a,b)) \mbox{ for } a,b\in A\},
$$
is a monoid with respect to composition. It is easy to check that this monoid is Zariski closed in $\LL(\vect(A))$, therefore it is an affine algebraic monoid.
 
 It is shown in \cite{GP} that any affine algebraic group can be realized as the group of automorphisms of some finite-dimensional simple algebra. This note aims to obtain  a similar realization of an 
 arbitrary affine algebraic monoid $M$ as the endomorphisms' monoid of a finite-dimensional algebra $A$. In this case two differences 
 occur. First, 
 we cannot assume that $A$ is simple, since the kernel of any endomorphism is an ideal of $A$. Second, the monoid $\End(A)$ contains 
  zero $\z\in\End(A)$, $\z(a)=0$ for any $a\in A$, while 
   $M$ does not necessarily contain 
   zero. Under these circumstances we obtain the following result.

\begin{theo}\label{main}
For any affine algebraic monoid $M$ there exists a finite-dimensional algebra $A$ such that
$\End(A)\cong M\sqcup\{\z\}$, where $\{\z\}$ is an (isolated) component of the monoid $\End(A)$.
\end{theo}

Particularly, if $M$ is an affine algebraic group, then there exists an algebra $A$ such that $\Aut(A)\cong M$ (see \cite{GP}).
\begin{example}
Let us consider the monoid $M=\LL(V)$ for a finite-dimensional space $V$. Then we may take the algebra $A$ constructed in the following way.
First, let $e$ be a left identity of $A$ and
$$
\vect(A)\coloneq\bangle{e}\oplus V,
$$
where $\bangle{X}$ stands for the linear span of a set $X$.
Next, for any $v,w\in V$ put $\alpha(v,w)=0$, $\alpha(v,e)=\lambda v$, where $\lambda\in\K\setminus\{0,1\}$. Taking into account equations $\alpha(e,v)=v$ and $\alpha(e,e)=e$, we obtain the multiplication table for $A$.

Note that any endomorphism sends $e$ to $e$ or $0$, since these two are the only idempotents of $A$. This way, the reader will easily prove that $\End(A)\cong\LL(V)\sqcup\{\z\}$.
\end{example}

\begin{example}Assume $\chr\K\neq2$.
Consider a two-dimensional space $V$ over $\K$ with a basis $\{v_1,v_2\}$ and the exterior algebra $\Lambda(V)$ with a basis $\{1,v_1,v_2,v_1\wedge v_2\}$. Let us take a monoid $M\subset\LL(\vect(\Lambda(V)))$,
$$
M\coloneq
\left\{\left.\left(\begin{array}{cccc}
1 & 0& 0& 0\\
0& b_{11} & b_{12} & 0\\
0& b_{21} & b_{22} & 0\\
0& c_1 & c_2 & d
\end{array}\right)\right|d=\det\left(\begin{array}{cc}b_{11} & b_{12}\\
b_{21} & b_{22}\end{array}\right), b_{ij}, c_i\in\K\right\}.
$$
One may prove that $M$ acts on $\Lambda(V)$ by endomorphisms.
Moreover, $\End(\Lambda(V))=M\sqcup\{\z\}$. Generally, a similar equation holds for the exterior algebra of an arbitrary space.
\end{example}

The proof of Theorem \ref{main} consists of two steps. First, for every finite-dimensional space $U$ and its subspace $S$ we construct a finite-dimensional algebra $A$ such that $\End(A)$ is isomorphic to $\LL(U)_S\sqcup\{\z\}$, where $\LL(U)_S$ is the normalizer of some vector subspace $S$ of a special $\LL(U)$-module.
Second, arbitrary affine algebraic monoid $M$ is represented as $\LL(U)_S$ for appropriate $U$ and $S$.
Overall, we follow the scheme of the proof in \cite{GP}, but the ideas of the first step are significantly changed.

\section{Some special algebras}
In this section we define and study some finite-dimensional algebras to be used hereafter.
\subsection{Algebra $A(V,S)$}
Let $V$ be a nonzero finite-dimensional vector space.
Denote by $\T(V)$ the tensor algebra of $V$
    and by $\T(V)_+$ its maximal homogeneous ideal
    \begin{equation}
    \T(V)_+\coloneq\boplus_{i\ge 1}V^{\otimes i},
    \end{equation}
    endowed with the natural $\LL(V)$-structure,
    \begin{equation}\label{L-act}
    g\cdot t_i\coloneq g^{\otimes i}(t_i),\quad g\in\LL(V),\: t_i\in V^{\otimes i}.
    \end{equation}
    Thus, $\LL(V)$ acts on $\T(V)_+$ faithfully by endomorphisms.
    Therefore we may  
    identify $\LL(V)$ with the corresponding submonoid of $\End(\T(V)_+)$.

Fix an integer $r>1$. For arbitrary subspace $S\subseteq V^{\otimes r}$ we define
 \begin{equation}\label{I(S)}
 I(S)\coloneq S\oplus(\boplus_{i>r}V^{\otimes i}).
 \end{equation}
 It is an ideal of $\T(V)_+$. Define $A(V,S)$ as the factor algebra modulo this ideal,
 \begin{equation}\label{A(V,S)}
 A(V,S)\coloneq \T(V)_+/I(S).
 \end{equation}
 Then
 \begin{equation}\label{vect(A)}
 \vect(A(V,S))=(\boplus_{i=1}^{r-1}V^{\otimes i})\oplus(V^{\otimes r}/S).
 \end{equation}

 We may consider $\LL(V)_S\coloneq\{\phi\in\LL(V)\;|\;\phi(S)\subseteq S\}\subset\LL(V)$.

 \begin{proposition}\label{End(A)}
 $\{\sigma\in\End(A(V,S))\;|\;\sigma(V)\subseteq V\}=\LL(V)_S$.
 \end{proposition}
 \begin{proof}
 By definition, elements of $A(V,S)$ are equivalence classes $x+I(S)$, $x\in \T(V)_+$.
Let us prove the inclusion $\subseteq$. Consider $\sigma\in\End(A(V,S))$ such that $\sigma(V)\subseteq V$.
  Then $\sigma$-action coincides with the action of $\widetilde{\sigma}\coloneq\sigma|_V\in\LL(V)$ on $A(V,S)$ in accordance with (\ref{L-act}), since the algebra $A(V,S)$ is generated by $V$. The $\sigma$-action preserves zero of $A(V,S)$, hence $\widetilde{\sigma}(I(S))\subseteq I(S)$ and $\sigma\in\LL(V)_S$.

 Now we prove the inverse inclusion. For arbitrary subsets $X,Y\subset \T(V)$ define $X\otimes Y\coloneq\{x\otimes y\;|\; x\in X,y\in Y\}\subset\T(V)$. Let $\sigma\in\LL(V)_S$. Then ${\sigma}((x+I(S))\otimes(y+I(S)))\subseteq {\sigma}(x\otimes y)+I(S)={\sigma}(x)\otimes {\sigma}(y)+I(S)$ by definition of the $\LL(V)$-action on $\T(V)_+$. Hence $\sigma\in\End(A(V,S))$.
 \end{proof}

\subsection{Algebra $D(P,U,S,\mathbf{\gamma})$}
\begin{lemma}\label{oplus-lemma}
Let $A$ be an algebra with a left identity $e\in A$ such that $\vect(A)=\bangle{e}\oplus A_1\oplus\cdots\oplus A_r$, where
  $A_i$ is the 
  eigenspace with an eigenvalue $\alpha_i\neq0,1$ for 
   the operator of right multiplication of $A$ by $e$. Assume that 0 and $e$ are the only idempotents in~$A$. Then
 \begin{enumerate}
 \item $e$ is the unique left identity in $A$;
 \item if $\sigma\in\End(A)$, then either $\sigma(e)=e$ and $\sigma(A_i)\subseteq A_i$ for any $i$, or $\sigma=\z$.
 \end{enumerate}
\end{lemma}
\begin{proof}
(i) The left identity is a nonzero idempotent. Hence it is unique.

(ii) Since the image of an idempotent is an idempotent, $\sigma(e)=0$ or $\sigma(e)=e$.
If $\sigma(e)=0$, then $\sigma(a)=\sigma(ea)=\sigma(e)\sigma(a)=0$, i.e. $\sigma=\z$.
Now assume that $\sigma(e)=e$.
Then $\sigma(A_i)$ is the eigenspace with an eigenvalue $\alpha_i\neq0,1$ for the operator of right multiplication by $e$. Hence $\sigma(A_i)\subseteq A_i$.
\end{proof}

Let $P$ be a two-dimensional vector space with a basis $\{p_1,p_2\}$, $U$ be a nonzero finite-dimensional space, and
 \begin{equation}
    V\coloneq P\oplus U.
 \end{equation}
Fix an integer $r>1$ as well as 
\begin{list}{(\roman{countroman})}{\usecounter{countroman}}
\item a subspace $S\subset V^{\otimes r}$;
\item a sequence $\mathbf{\gamma}=(\gamma_1,\ldots,\gamma_6)\in(\K\setminus\{0,1\})^6,\;\gamma_i\neq\gamma_j$ for $i\neq j$.
\end{list}

Define an algebra $D(P,U,S,\mathbf{\gamma})$ in the following way.
First, $A(V,S)$ is the subalgebra of $D(P,U,S,\mathbf{\gamma})$ and elements $b,c,d,e\in D(P,U,S,\mathbf{\gamma})$ are such that
 \begin{equation}\label{vect(D)}
    \vect(D(P,U,S,\mathbf{\gamma}))=\bangle{e}\oplus\bangle{b}\oplus\bangle{c}\oplus\bangle{d}\oplus\vect(A(V,S)).
 \end{equation}

 Second, the following conditions hold:
\begin{list}{(D\arabic{count})}{\usecounter{count}}
\item $e$ is the left identity of $D(P,U,S,\mathbf{\gamma})$;
\item $\bangle{b}$, $\bangle{c}$, $\bangle{d}$ as well as $P,U\subset V=V^{\otimes1}\subset A(V,S)$ and $(\boplus_{i=2}^{r-1} V^{\otimes i})\oplus(V^{\otimes r}/S)\subset A(V,S)$ are the eigenspaces with the eigenvalues $\gamma_1,\ldots,\gamma_6$ respectively of the operator of right multiplication by $e$;
\item The multiplication table for $b,c,d$ is\label{D3} 
\begin{equation}\begin{array}{lll}
b\cdot b\coloneq0    ,\;& b\cdot c\coloneq c+\gamma_{bc}b  ,\;& b\cdot d\coloneq0  ,\\
c\cdot b\coloneq -c   ,\; & c\cdot c\coloneq b ,\;& c\cdot d\coloneq e ,\\
d\cdot b\coloneq p_1    ,\;& d\cdot c\coloneq d  ,\;& d\cdot d\coloneq p_2  ,
\end{array}\end{equation}
where $\gamma_{bc}=\frac{\gamma_2-\gamma_1}{\gamma_2-\gamma_3}$;
\smallskip
\item $\bangle{b,c,d}\cdot A(V,S)=A(V,S)\cdot \bangle{b,c,d}=0$.
\end{list}

Define the action of $g\in\LL(V)_S$ on $\vect(D(P,U,S,\mathbf{\gamma}))$ as follows: $g|_{\bangle{b}}=g|_{\bangle{c}}=g|_{\bangle{d}}=g|_{\bangle{e}}=\id$,
$g|_V$ is the natural $\LL(V)$-action on $V$, and on other summands of $A(V,S)$ it is defined by (\ref{L-act}). By Proposition \ref{End(A)} we may identify $\LL(V)_S$ with the corresponding submonoid of $\LL(\vect(D(P,U,S,\mathbf{\gamma})))$.
Further, we may consider an embedding $\LL(U)\embed\LL(V)$, $h\mapsto \id|_P\oplus h$. Thus, $\LL(U)_S\subseteq\LL(V)_S$, and we obtain
$\LL(U)_S$-action on $\vect(D(P,U,S,\mathbf{\gamma}))$.

\begin{proposition}\label{End(D)}
We have 
$$\End(D(P,U,S,\mathbf{\gamma}))=\LL(U)_S\sqcup\{\z\},$$ where $\{\z\}$ is an (isolated) component of the monoid $\End(D(P,U,S,\mathbf{\gamma}))$.
\end{proposition}
\begin{proof}
First of all, we show that 0 and $e$ are the only idempotents of $D(P,U,S,\mathbf{\gamma})$.
Indeed, let $\eps=\lambda_e e+\lambda_b b +\lambda_c c +\lambda_d d +a$, where $a\in A(V,S)$. Then
\begin{multline}
\eps^2=(\lambda_e^2+\lambda_c\lambda_d) e+(\lambda_b\lambda_e(1+\gamma_1)+\lambda_c^2+\lambda_b \lambda_c\gamma_{bc}) b +\\
+\lambda_c\lambda_e(1+\gamma_2) c +((1+\gamma_3)\lambda_d\lambda_e+\lambda_d\lambda_c) d +a'=\\
=\lambda_1 e+\lambda_2 d +\lambda_3 c +\lambda_4 b +a,\; \mbox{where}\; a,a'\in A(V,S).
\end{multline}
Hence
\begin{gather}
\lambda_e=\lambda_e^2+\lambda_c\lambda_d,\label{lambda-e}\\
\lambda_b=\lambda_b\lambda_e(1+\gamma_1)+\lambda_c^2+\lambda_b \lambda_c\gamma_{bc},\label{lambda-b}\\
\lambda_c=\lambda_c\lambda_e(1+\gamma_2),\label{lambda-c}\\
\lambda_d=\lambda_d\lambda_e(1+\gamma_3)+\lambda_c\lambda_d.\label{lambda-d}
\end{gather}

Assume $\lambda_c\neq0$. By (\ref{lambda-c}) $1+\gamma_2\neq0, \lambda_e=\frac{1}{1+\gamma_2}$ and $\lambda_c\lambda_d=\lambda_e-\lambda_e^2\neq0$, so $\lambda_d\neq0$. Hence equation (\ref{lambda-d}) implies $\lambda_c=1-\lambda_e(1+\gamma_3)=\frac{\gamma_2-\gamma_3}{1+\gamma_2}$. Finally, by (\ref{lambda-b}) we have $\lambda_c^2=\lambda_b(1-\lambda_e(1+\gamma_1)-\lambda_c\gamma_{bc})=\lambda_b(\frac{\gamma_2-\gamma_1}{1+\gamma_2}-\frac{\gamma_2-\gamma_3}{1+\gamma_2} \cdot\frac{\gamma_2-\gamma_1}{\gamma_2-\gamma_3})=0$. From this contradiction we deduce $\lambda_c=0$.

   Moreover, $\lambda_e=0$ or $\lambda_e=1$ by (\ref{lambda-e}). If $\lambda_e=0$, then $\lambda_b=\lambda_d=0$, $\eps=a\in A(V,S)$ and $\eps=0$, since zero is the only idempotent of $A(V,S)$. Now assume $\lambda_e=1$. From equations (\ref{lambda-b}) and (\ref{lambda-d}) accordingly follow
  $\lambda_b=0$ and $\lambda_d=0$. Thus, $\eps=e+a$, $a\in A(V,S)$.

Let $a=a_P+a_U+a_{\Sigma}$, where $a_P\in P, a_U\in U, a_{\Sigma}\in (\boplus_{i=2}^{r-1} V^{\otimes i})\oplus(V^{\otimes r}/S)$. Then
\begin{equation}
\eps^2=e+(1+\gamma_4)a_P+(1+\gamma_5)a_U+a'_{\Sigma}=e+a_P+a_U+a_{\Sigma},
\end{equation}
where $a'_{\Sigma}\in (\boplus_{i=2}^{r-1} V^{\otimes i})\oplus(V^{\otimes r}/S)$.
Hence $a_U=a_P=0$. Assume $a_{\Sigma}\neq0$, then we may write $a_{\Sigma}=a_k+\ldots+a_r$, $a_k\neq0$, where $a_i\in V^{\otimes i}$ for $i<r$ and $a_r\in V^{\otimes r}/S$.
This way,
\begin{equation}
(e+a_k+\ldots+a_r)^2=
e+(1+\gamma_6)a_k+a''=
e+a_k+\ldots+a_r,
\end{equation}
where $a''\in(\boplus_{i=k+1}^{r-1} V^{\otimes i})\oplus(V^{\otimes r}/S)$ for $k<r$ and $a''=0$ for $k=r$.
It implies $a_k=0$. Contradiction. Hence $a_{\Sigma}=0$ and $\eps=e$.

Thus, $D(P,U,S,\mathbf{\gamma})$ contains no idempotents different from 0 and $e$. Let $\sigma\in\End(D(P,U,S,\mathbf{\gamma}))\setminus\{\z\}$.
By Lemma \ref{oplus-lemma} $\sigma(e)=e$ and $\bangle{b}$, $\bangle{c}$, $\bangle{d}$, $P$, $U$, $A(V,S)$ are $\sigma$-invariant.
Let $\sigma(b)=\delta_b b$, $\sigma(c)=\delta_c c$, $\sigma(d)=\delta_d d$.
The equations $cd=e,\:dc=d,\:cb=-c$ imply $\delta_c\delta_d=1,\:\delta_c\delta_d=\delta_d,\:\delta_b\delta_c=\delta_c$.
One may check that $\delta_b=\delta_c=\delta_d=1$. Finally, the equations $db=p_1,\:dd=p_2$ imply $\sigma|_P=\id_P$.

 Since $V$ and $A(V,S)$ are $\sigma$-invariant, $\sigma|_{A(V,S)}\in\LL(V)_S$ by Proposition \ref{End(A)}. Taking into account $\sigma|_P=\id_P$
and $\sigma(U)\subseteq U$, we obtain $\sigma\in\LL(U)_S$.
\end{proof}

\section{Affine monoids as the normalizers of linear subspaces}
\begin{proposition}\label{normalization}
Let $M$ be an affine algebraic monoid.
There is a finite-dimensional vector space $U$ and an integer $r>1$ such that the following holds.
Let $P$ be a two-dimensional vector space with a trivial $\LL(U)$-action.
Then the $\LL(U)$-module $(P\oplus U)^{\otimes r}$ contains a linear subspace $S$ such that $\LL(U)_S\cong M$.
\end{proposition}
\begin{proof}
Since there exists a closed embedding $M\embed\LL(U)$ for some finite-dimensional space $U$, 
we may suppose $M\subseteq\LL(U)$.
Consider the action of $\LL(U)$ on itself by left multiplication.
Additionally, consider the $\LL(U)$-action on the algebra $\K[\LL(U)]$ of regular 
functions on $\LL(U)$,
\begin{equation}
(g\cdot f)(u)\coloneq f(ug),\quad g,u\in\LL(U), f\in\K[\LL(U)].
\end{equation}

  Denote $d\coloneq\dim U$.
  Note that the $\LL(U)$-modules $\K[\LL(U)]$ and $\Sym(U^{\oplus d})$ are isomorphic. To prove this, it suffices to associate a linear function on $\LL(U)$ to every vector $(u_1,\ldots,u_d)\in U^{\oplus d}$, since $\K[\LL(U)]=\Sym(\LL(U)^*)$. Identify $U$ with $\K^d$, $\LL(U)$ with $\Mat_{d\times d}(\K)$; let $A$ be in $\LL(U)$, $B$ be a matrix with columns $u_1,\ldots,u_d$. Set $l_{u_1,\ldots,u_d}(A)\coloneq\tr AB$. Then $(g\cdot l_{u_1,\ldots,u_d})(A)=\tr AgB=l_{gu_1,\ldots,gu_d}(A)$, i.e. we have an $\LL(U)$-equivariant isomorphism.

By definition of symmetric algebra fix a natural epimorphism
 \begin{equation}
 \xi\colon\T(U^{\oplus d})\to\Sym(U^{\oplus d})\cong\K[\LL(U)].
 \end{equation}

There is a finite-dimensional subspace $W\subset \K[\LL(U)]$ such that
 \begin{equation}\label{LL(U)_W}
 \LL(U)_W=M.
 \end{equation}

In order to prove this, one may show that a linear span of an $\LL(U)$-`orbit' of an arbitrary function $f\in\K[\LL(U)]$ is finite-dimensional.
 Indeed, since the $\LL(U)$-action is a morphism, 
  $(g\cdot f)(u)=f(ug)\in\K[\LL(U)\times \LL(U)]=\K[\LL(U)]\otimes \K[\LL(U)]$, where $u,g\in\LL(U)$,
 there are functions $F_j,H_j\in\K[\LL(U)]$ such that
 \begin{equation}
     (g\cdot f)(u)=\sum_{j=1}^n F_j(u)H_j(g).
 \end{equation}
Therefore, the $\LL(U)$-`orbit' of the function $f$ is contained in the finite-dimensional subspace $\bangle{F_1,\ldots,F_n}$.

Let $I(M)=(f_1,\ldots,f_t)\lhd\K[\LL(U)]$ be the ideal of functions vanishing on $M$.
Summing the linear spans of $\LL(U)$-`orbits' of the functions $f_i$ we obtain a finite-dimensional $\LL(U)$-invariant subspace $V\subset\K[\LL(U)]$. Define $W=I(M)\cap V$. First, it contains $f_1,\ldots,f_t$. Second, it is $M$-invariant, since the ideal $I(M)$ is $M$-invariant.
 Obviously, $g\in M$ implies $g\cdot W\subseteq W$.
  On the other hand, suppose that $g\cdot W\subseteq W$, where $g\in\LL(U)$.
   Then $f_i(g)=(g\cdot f_i)(E)=0$ for $i=1,\ldots,t$, where $E$ is the identity of $\LL(U)$ and is automatically contained in $M$.
  Therefore, $g\in M$. This proves (\ref{LL(U)_W}).

 Further, since the space $W$ is finite-dimensional, there is an integer $h\in\ZZ_+$ such that
 \begin{equation}
 W\subseteq\xi(\boplus_{i\le h}(U^{\oplus d})^{\otimes i}).
 \end{equation}
 Define $W'\coloneq\xi^{-1}(W)\cap(\boplus_{i\le h}(U^{\oplus d})^{\otimes i})$. The $\LL(U)$-equivariance of $\xi$ implies
 \begin{equation}\label{LL(U)_W'}
 \LL(U)_{W'}=\LL(U)_W.
 \end{equation}

 Fix a basis $\{p_1,p_2\}$ of the space $P$.
 There exists an embedding of $\LL(U)$-modules
 \begin{equation}
 \iota\colon\T(U^{\oplus d})\embed\T(\bangle{p_1}\oplus U).
 \end{equation}
 Indeed, let $U_i$ be the ith summand of $U^{\oplus d}$. Consider an arbitrary basis $\{f_{i j}\;|\;j=1,\ldots,d\}$
 of $U_i$ and define an embedding as follows,
 \begin{equation}\label{i1}
 \iota(f_{i_1 j_1}\otimes\ldots\otimes f_{i_t j_t})\coloneq p_1^{\otimes i_1}\otimes f'_{i_1 j_1}\otimes\ldots
 \otimes p_1^{\otimes i_t}\otimes f'_{i_t j_t},
 \end{equation}
 where $f'_{i j}$ is the image of $f_{i j}$ under the identity isomorphism $U_i\to U$. It is easy to check that the embedding
 $\iota\colon\T(U^{\oplus d})\to\T(\bangle{p_1}\oplus U)$ defined on the basis of $\T(U^{\oplus d})_+$ by formula (\ref{i1})
 and sending 1 to 1 is the one required. 

 Now we may consider a space $W''\coloneq\iota(W')$,
 \begin{equation}\label{LL(U)_W''}
 \LL(U)_{W''}=\LL(U)_{W'}.
 \end{equation}
 Since $W''$ is finite-dimensional, there exists an integer $b\in\NN$ such that
 \begin{equation}
 W''\subseteq\boplus_{i\le b}(\bangle{p_1}\oplus U)^{\otimes i}.
 \end{equation}
 Take $r\ge b$ such that 
  $r>1$ and consider a linear mapping
 \begin{equation}\label{i2}
 \iota_r\colon\boplus_{i\le b}(\bangle{p_1}\oplus U)^{\otimes i}\to(P\oplus U)^{\otimes r},\quad
 f_i\mapsto p_2^{\otimes (r-i)}\otimes f_i, f_i\in(\bangle{p_1}\oplus U)^{\otimes i}.
 \end{equation}
 Obviously, $\iota_r$ is an embedding of $\LL(U)$-modules. Define $S=\iota_r(W'')$. Then
 \begin{equation}\label{LL(U)_S}
 \LL(U)_S=\LL(U)_{W''}.
 \end{equation}
 Now the claim follows from equations (\ref{LL(U)_W}), (\ref{LL(U)_W'}), (\ref{LL(U)_W''}) and (\ref{LL(U)_S}).
\end{proof}

\subsection*{Proof of Theorem \ref{main}}
Let $M$ be an arbitrary affine algebraic monoid, $U,b,r,P,S$ be as in Proposition \ref{normalization}.
Fix some set $\mathbf{\gamma}\in(\K\setminus\{0,1\})^6$ such that $\gamma_i\neq\gamma_j$ for $i\neq j$,
and consider the algebra $D(P,U,S,\mathbf{\gamma})$.
It follows from Proposition \ref{normalization} and Proposition \ref{End(D)} that $\End(D(P,U,S,\mathbf{\gamma}))\cong M\sqcup\{\z\}$.

\section*{Acknowledgements}
The author expresses his sincere thanks to I.V.~Arzhantsev for posing of the problem and useful discussions.


\begin{thebibliography}{}%
%
\bibitem{A} I.V.~Arzhantsev, \emph{Affine embeddings of homogeneous spaces},
in {\it ``Surveys in Geometry and Number Theory''}, N.~Young (Editor),
LMS Lecture Notes Series {\bf 338}, Cambridge Univ. Press, Cambridge, 2007, 1--51
%
\bibitem{GP} N.L.~Gordeev and V.L.~Popov, \emph{Automorphism groups of finite dimensional simple algebras},
 Annals of Mathematics {\bf 158} (2003), 1041--1065
%
\bibitem{P} M.S.~Putcha, \emph{Linear algebraic monoids}, LMS Lecture Notes Series \textbf{133}, Cambridge Univ. Press, Cambridge, 1988
%
\bibitem{R} L.~Renner, \emph{Linear algebraic monoids}, Encyclopaedia of Mathematical Sciences \textbf{134}, Springer-Verlag, Berlin Heidelberg, 2005
%
\bibitem{Rit} A.~Rittatore, \emph{Algebraic monoids and affine embeddings}, Transform. Groups \textbf{3} (1998) no. 4, 375--396
%
\bibitem{V} E.B.~Vinberg, \emph{On reductive algebraic semigroups}, Amer. Soc. Transl. (2) \textbf{169} (1995), 145--182
%
\end{thebibliography}
\end{document}